\documentclass{article}
\usepackage{maa-monthly}
\usepackage[final]
{showkeys}
\usepackage{graphicx}
\usepackage{enumerate}
\usepackage{multirow}
\usepackage[skip=40pt,font=large]{caption}
\usepackage{framed}

\usepackage[framemethod=tikz]{mdframed}

\usepackage{trimclip}

\usepackage{needspace}

\usepackage{esint}

\usepackage{hyperref}

\theoremstyle{theorem}
\newtheorem{theorem}{Theorem}

\newtheorem{corollary}[theorem]
{Corollary}

{Claim}

\theoremstyle{definition}

\newtheorem{remark}[theorem]
{Remark}

\theoremstyle{definition} 
\newtheorem*{remark*}{Remark}

\renewcommand{\skip}[1]{}

\newcommand{\R}{\mathbb{R}}

\newcommand{\de}{\delta}

\newcommand\la{\lambda}

\newcommand\co{\operatorname{conv}}
\newcommand\ol{\overline}

\makeatletter
\DeclareRobustCommand{\lasymp}{\lg@asymp{<}}
\DeclareRobustCommand{\gasymp}{\lg@asymp{>}}

\newcommand{\under@asymp}[1]{\clipbox{0pt 0pt 0pt {0.5\height}}{$\m@th#1\asymp$}}
\newcommand{\lg@asymp}[1]{\mathrel{\mathpalette\lg@asymp@{#1}}}
\newcommand{\lg@asymp@}[2]{%
  \vcenter{%
    \offinterlineskip
    \m@th
    \ialign{%
      \hfil##\hfil\cr
      $#1#2$\cr
      \under@asymp{#1}\cr
    }%
  }%
}
\makeatother

\makeatletter

\long\def\@maketablecaption#1#2{%
    \rule[.1\baselineskip]{0pt}{\baselineskip}
    \small\textbf{#1.}\enspace #2\strut
    \par
  \vskip2pt}
\makeatother  

\begin{document}

\title{Does the Convex Order Between the Distributions of Linear Functionals Imply the Convex Order Between the Probability Distributions Over $\R^d$? 
}
\markright{
\quad\quad\quad\quad\quad\ \  Convex order over $\mathbb R^d$ via projections}
\author{Iosif Pinelis \\ 
Michigan Technological University \\ 
Houghton, Michigan USA \\ 
Email: ipinelis@mtu.edu 
}


%
%
%
%
%

\maketitle

\begin{abstract} 
It is shown that the convex order between the distributions of linear functionals does not imply the convex order between the probability distributions over $\mathbb R^d$ if $d\ge2$. This stands in contrast with the well-known fact that any probability distribution in $\R^d$, for any $d\ge1$, is  determined by the corresponding distributions of linear functionals. By duality, it follows that, for any $d\ge2$, not all convex functions from $\R^d$ to $\R$ can be represented as the limits of sums $\sum_{i=1}^k g_i\circ \ell_i$ of convex functions $g_i$ of linear functionals $\ell_i$ on $\mathbb R^d$.  
\end{abstract}










Let $\mu$ and $\nu$ be probability measures over $\R^d$ with 
\begin{equation}\label{eq:finite}
\text{$\int_{\R^d}\|x\|\,\mu(dx)<\infty$\quad and\quad $\int_{\R^d}\|x\|\,\nu(dx)<\infty$,}	
\end{equation}
where $\|\cdot\|$ denotes the Euclidean norm. 
It is said that $\mu$ is dominated by $\nu$ in the convex order sense and written $\mu\preceq\nu$ if 
\begin{equation}\label{eq:}
	\int_{\R^d}f\,d\mu\le\int_{\R^d}f\,d\nu
\end{equation}
for all convex functions $f\colon\R^d\to\R$. 
Note that, given \eqref{eq:finite}, the integrals in \eqref{eq:} always exist for any convex function $f\colon\R^d\to\R$ but may take the value $\infty$---because for any such $f$ there is an affine function $a\colon\R^d\to\R$ such that $f\ge a$. 

For any affine function $a\colon\R^d\to\R$, both $a$ and $-a$ are convex functions. So, if $\mu\preceq\nu$, then the barycenters $x_\mu:=\int_{\R^d}x\,\mu(dx)$ and $x_\nu:=\int_{\R^d}x\,\nu(dx)$ of $\mu$ and $\nu$ are the same. 

Informally, the convex order relation $\mu\preceq\nu$ means that, while the barycenters $x_\mu$ and $x_\nu$ of the distributions $\mu$ and $\nu$ are the same, the distribution $\nu$ is more spread out than $\mu$. For instance, one has $\de_y\prec\nu$ for any $y\in\R^d$ and any probability measure $\nu$ with barycenter $y$, where $\de_y$ is the Dirac probability measure supported on the singleton set $\{y\}$. 

The convex order has been widely studied and/or used; see, e.g., Jovan Karamata \cite{karamata32}, David Blackwell \cite{blackwell53}, Paul-Andr\'e Meyer \cite{meyer}, and Robert Ralph Phelps \cite{phelps}. For stochastic orders in general, see Moshe Shaked and J. George Shanthikumar \cite{shaked-shanti}. 
The notion of the decreasing concave order---which is a reverse of the increasing convex order and is also known as the second-order stochastic dominance---was introduced by Michael Rothschild and Joseph Eugene Stiglitz \cite{rothschild-stiglitz} and has been very widely used in economics literature afterwards. 

For any probability measure $\la$ over $\R^d$ and any $v\in\R^d$, let $\la_v$ denote the \break  ``$v$-projection'' of $\la$ -- that is, the pushforward of $\la$ under the linear map   
\begin{equation}\label{eq:proj}
\R^d\ni x\mapsto p_v(x):=v\cdot x\in\R, 	
\end{equation}
where $\cdot$ stands for the dot product. 

It is well known that any probability measure $\la$ over $\R^d$ is determined by the family $(\la_v)_{v\in\R^d}$ of its one-dimensional projections; see, e.g., \cite[Theorem 26.2]{billingsley95}. 
Using this result, one can establish its extension, due to Cram\'er and Wold \cite{cramer-wold}, stating that the weak convergence of multivariate distributions is determined by the weak convergence of their one-dimensional projections; see also, e.g., \cite[Theorem 29.4]{billingsley95}. 

One may then ask whether the convex order is similarly determined by the one-dimensional projections---that is, whether 
\begin{equation}\label{eq:iff}
	\mu\preceq\nu \overset{\text{(?)}}\iff \big(\mu_v\preceq\nu_v \text{ for all }v\in\R^d\big). 
\end{equation}
The implication $\implies$ in \eqref{eq:iff} is obvious, because the composition $g\circ p_v$ is convex for all convex functions $g\colon\R\to\R$ and all $v\in\R^d$; of course, for the same reason, an appropriate version of this implication holds for any real topological vector space in place of $\R^d$. 

So, the question is only about the implication $\impliedby$ in \eqref{eq:iff}. 
This implication is trivially true if $d=1$.  
However, we will show that the answer 
to the question 
is negative 
if $d\ge2$. 
In fact, a counterexample to the implication $\impliedby$ in \eqref{eq:iff} will be given explicitly by two probability measures $\mu$ and $\nu$ that have finite support sets and take rational values. 
Of course, it follows that the answer to the corresponding question for any real topological vector space of dimension $d\ge2$ in place of $\R^d$ will also be negative. 
 
%
%
%
%
%
%
%

Indeed, already for $d=2$, let 
\begin{equation}\label{0}
	\nu=\frac16\,(3\de_{(x+y+z)/3}+\de_x+\de_y+\de_z),\quad 
\mu=\frac13\,(\de_{(x+y)/2}+\de_{(y+z)/2}+\de_{(x+z)/2}) 
\end{equation}
for some $x,y,z$ in $\R^2$. 
Then $\mu\preceq\nu$ means that 
\begin{equation}\label{1}
	3f(\tfrac{x+y+z}3)+f(x)+f(y)+f(z)
\ge2f(\tfrac{x+y}2)+2f(\tfrac{y+z}2)+2f(\tfrac{x+z}2) 
\end{equation}
for all convex $f\colon\R^2\to\R$. 
Similarly, for each $v\in\R^2$, the relation $\mu_v\preceq\nu_v$ means that \eqref{1} holds for all functions $f$ of the form $g\circ p_v$, where $g$ is any convex function from $\R$ to $\R$ and $p_v$ is as defined in \eqref{eq:proj}. 
So, for each $v\in\R^2$, the relation $\mu_v\preceq\nu_v$ 
means that  
\begin{equation}\label{2}
	3g(\tfrac{r+s+t}3)+g(r)+g(s)+g(t)
\ge2g(\tfrac{r+s}2)+2g(\tfrac{s+t}2)+2g(\tfrac{r+t}2) 
\end{equation}
for $r=p_v(x)$, $s=p_v(y)$, and $t=p_v(z)$ and all convex functions $g\colon\R\to\R$.   

In turn, inequality \eqref{2} does hold for all such convex functions $g$ and all real $r,s,t$---being the simplest instance of an 
inequality due to Tiberiu Popoviciu \cite{popoviciu65}; see also e.g.\ \cite[Theorem 1.1.8]{nicul-perss} and \cite [p.\ 74]{fechner14}. Inequality \eqref{2} is also a special case (with $p=q=r=1/3$) of \cite[inequality (6.2)]{pecaric-proschan-tong}.

It is easy to prove \eqref{2} directly as well. To do this, let us first recall the definition and a basic characterization 
of majorization. For a vector $x=(x_1,\dots,x_n)$ in $\R^n$, let $
(x_{[1]},\dots,x_{[n]})$ be the nonincreasing rearrangement of the $x_i$'s, so that 
$x_{[1]}\ge\dots\ge x_{[n]}$. Let us then say that a vector $x=(x_1,\dots,x_n)\in\R^n$ is majorized by a vector $y=(y_1,\dots,y_n)\in\R^n$ (and write $x\prec y$) if $\sum_{i=1}^k x_{[i]}\le\sum_{i=1}^k y_{[i]}$ for all $k=1,\dots,n-1$ and $\sum_{i=1}^n x_{[i]}=\sum_{i=1}^n y_{[i]}$---see e.g.\ \cite[Definition 1.A.1]{marsh-ol11}. 

A most important characterization of majorization is as follows: For vectors $x$ and $y$ as above, one has 
$x\prec y$ if and only if $\sum_{i=1}^n g(x_i)\le\sum_{i=1}^n g(y_i)$ for all convex functions $g\colon\R\to\R$ (cf., e.g.,  \cite[Proposition~4.B.1]{marsh-ol11}, where it is additionally required that $g$ be continuous; however, any convex function from $\R$ to $\R$ is automatically continuous---see, e.g.,  \cite[Corollary~10.1.1]{rocka}).  

%

%

So, to prove \eqref{2}, it suffices to check that  
\begin{equation}\label{eq:R<L}
(\tfrac{r+s}2,\tfrac{r+s}2,\tfrac{s+t}2,\tfrac{s+t}2,\tfrac{r+t}2,\tfrac{r+t}2)
\prec(\tfrac{r+s+t}3,\tfrac{r+s+t}3,\tfrac{r+s+t}3,r,s,t)
\end{equation} 
for any real $r,s,t$. 
To do this, first note that here, in view of permutation symmetry, without loss of generality $r\le s\le t$. Also, in view of the 
reflection symmetry 
$u\leftrightarrow-u$, without loss of generality $s\le\frac{r+t}2$ and hence 
$$r\le\tfrac{r+s}2\le s\le\tfrac{r+s+t}3\le\tfrac{r+t}2\le\tfrac{s+t}2\le t.$$ 
Now the majorization \eqref{eq:R<L} is easy to check just by the definition. 

So, for any $x,y,z$ in $\R^2$, we do have $\mu_v\preceq\nu_v$ for all $v\in\R^2$. 

However, \eqref{1} fails to hold if, e.g., 
\begin{equation}\label{eq:f}
\text{$f(w)=\max(0,\xi_1,\xi_2)$ for all $w=(\xi_1,\xi_2)\in\R^2$}	
\end{equation}
and 
\begin{equation}\label{eq:xyz}
	x=(0,-1),\quad y=(-1,0),\quad z=(2,2). 
\end{equation}  
So, $\mu\not\preceq\nu$. 
$\qed$ 

\begin{remark}\label{rem:1}
One may note that in this counterexample $\nu$ is the uniform distribution on the multiset consisting of the vertices $x,y,z$ of a triangle (each vertex taken with multiplicity $1$) and the barycenter $\tfrac{x+y+z}3$ (taken with multiplicity $3$), whereas $\mu$ is the uniform distribution on the set of the midpoints of the sides of that triangle. 



\end{remark}

%

\begin{remark}\label{rem:2}
Since $\R^2$ can be linearly embedded into $\R^d$ for any natural $d\ge2$, clearly \eqref{eq:iff} fails to hold for any such $d$. 
\end{remark}

Simple duality arguments lead to the following corollary. 

\begin{corollary}\label{cor:}
For any integer $d\ge2$, there is a convex function $f\colon\R^d\to\R$ that is not in the closed (say, with respect to the topology of pointwise convergence) convex hull $\ol{\co F}$ of the set $F$ of all functions of the form $g\circ p_v$, where $g\colon\R\to\R$ is a convex function and $v\in\R^d$.   
\end{corollary}

Indeed, as in Remark~\ref{rem:2}, here without loss of generality $d=2$. 
Let then $\mu$ and $\nu$ be as in \eqref{0} with $x,y,z$ as in \eqref{eq:xyz}. Then, as was shown above, \eqref{eq:} holds for all $f\in F$ and hence for all $f\in\ol{\co F}$. It was also shown that \eqref{eq:} does not hold for the convex function $f\colon\R^2\to\R$ defined by \eqref{eq:f}. So, the latter function $f$ is not in $\ol{\co F}$. $\qed$

\medskip
\hrule
\medskip

This note is related to the previous one \cite{norm-decomp_publ}, sharing with it references \cite{nicul-perss,fechner14}.



\medskip
\hrule
\medskip

{\bf Disclosure of interest } No conflict of interest to declare. 

{\bf Disclosure of funding } No funding was received for this work.

%
\bibliography{C:/Users/ipinelis/Documents/pCloudSync/mtu_pCloud_02-02-17/bib_files/citations04-02-21}
%
%
\bibliographystyle{vancouver}

\nopagebreak
\nopagebreak

\end{document}